\documentclass[10pt,reqno]{amsart}
\usepackage{amsmath}
\usepackage{amssymb}
\usepackage{mathrsfs}
\usepackage{xcolor, enumitem, comment, todonotes}
\usepackage[all]{xy}
 \usepackage{url}
 \usepackage{amsthm}
 \usepackage{thmtools}
\usepackage{thm-restate}
\usepackage{hyperref,soul}

\usepackage[T1]{fontenc}

\theoremstyle{plain}
\newtheorem*{theorem*}{Theorem}

\newtheorem{theorem}{Theorem}[section]
\newtheorem{claim}[theorem]{Claim}

\newtheorem{lemma}[theorem]{Lemma}

\newtheorem{proposition}[theorem]{Proposition}

\theoremstyle{definition}
\newtheorem{definition}[theorem]{Definition}

\newtheorem{question}[theorem]{Question}

\theoremstyle{remark}
\newtheorem{remark}[theorem]{Remark}

\newtheorem{fact}[theorem]{Fact}

\newcount\skewfactor
\def\mathunderaccent#1#2 {\let\theaccent#1\skewfactor#2
\mathpalette\putaccentunder}
\def\putaccentunder#1#2{\oalign{$#1#2$\crcr\hidewidth
\vbox to.2ex{\hbox{$#1\skew\skewfactor\theaccent{}$}\vss}\hidewidth}}

\def\smallbox#1{\leavevmode\thinspace\hbox{\vrule\vtop{\vbox
   {\hrule\kern1pt\hbox{\vphantom{\tt/}\thinspace{\tt#1}\thinspace}}
   \kern1pt\hrule}\vrule}\thinspace}

\title{Two-cardinal Kurepa Hypotheses}
\author{Fanxin Wu}
\address{Department of Mathematical Sciences, Rutgers University, Piscataway, NJ 08854}
\email{fw173@rutgers.edu}

\date{\today}

\begin{document}

\begin{abstract}
We consider the two-cardinal Kurepa Hypothesis $\mathsf{KH}(\kappa,\lambda)$. We observe that if $\kappa\leq\lambda<\mu$ are infinite cardinals then $\lnot\mathsf{KH}(\kappa,\lambda)\land\mathsf{KH}(\kappa,\mu)\rightarrow\mathsf{KH}(\lambda^+,\mu)$, and show that in some sense this is the only $\mathsf{ZFC}$ constraint. The case of singular $\lambda$ and its relation to Chang's Conjecture and scales is discussed. We also extend an independence result about Kurepa and Aronszajn trees due to Cummings to the case of successors of singular cardinal.
\end{abstract}

\maketitle

\section{Introduction}

At the end of the first volume of \textit{Fundamenta Mathematicae} \cite{suslin1920} was a list of open problems in set theory and descriptive set theory. The third of them, attributed to Mikhail Suslin, was the following: it had been known since Cantor that the real line $\mathbb{R}$ is the unique dense linear order without endpoints that is complete and separable; can separability in the above characterization be weakened to the countable chain condition (ccc), i.e., any collection of disjoint nonempty open sets is countable? According to \cite{Kanamori_2011}, Suslin's problem would grow in significance partly because all other problems on the list were eventually solved.

A counterexample to Suslin's problem is called a Suslin line. \DJ uro Kurepa \cite{kurepa1936}, independently Edwin Miller and Wac\l aw Sierpi\'{n}ski, reformulated Suslin's problem in terms of tree: there exists a Suslin line if and only if there exists a Suslin tree, which is a tree of height $\omega_1$ that has no uncountable chain or antichain. In the 60s and 70s, Suslin's problem was finally shown to be independent of the axioms of $\mathsf{ZFC}$ set theory; the tree formulation proved valuable in the consistency proofs.

The tree formulation also suggested other natural problems related to trees. For example, in contrast to a Suslin tree which has no cofinal branch, a Kurepa tree is a tree of height $\omega_1$ with countable levels and at least $\aleph_2$ many cofinal branches. Combinatorial problems surrounding trees have inspired a large body of study in set theory.

This paper is mainly concerned with the two-cardinal version of Kurepa tree. Two-cardinal generalizations of many combinatorial principles have been considered \cite{donder-matet1993,shelah-shioya2006,cody-white2024}, where a structure indexed by $\kappa$ is replaced by one indexed by $\mathcal{P}_\kappa(\lambda)$. A prominent example is two-cardinal Aronszajn tree, which was first studied by Jech \cite{jech1973mess} (who also considered two-cardinal Suslin trees) and Magidor \cite{magidor1974}, and more recently by Weiss \cite{weiss2010subtle} and others. Two-cardinal Kurepa Hypotheses first appeared in \cite{Jensen-Kunen}, and have been studied by, e.g., Friedman and Golshani \cite{friedman-golshani2012independence}.

As with the one-cardinal case, existence of two-cardinal Kurepa trees is equivalent to existence of two-cardinal Kurepa families, which are slightly easier to define. A $(\kappa,\lambda)$-Kurepa family is a collection $\mathcal{F}$ of subsets of $\lambda$ such that $\forall x\in\mathcal{P}_\kappa(\lambda),\ |\{A\cap x:A\in\mathcal{F}\}|<\kappa$ and $|\mathcal{F}|\geq\lambda^+$. The two-cardinal Kurepa Hypothesis $\mathsf{KH}(\kappa,\lambda)$ says there exists a $(\kappa,\lambda)$-Kurepa family. Our first observation is the following simple $\mathsf{ZFC}$ restriction on combinations of different Kurepa Hypotheses: if $\kappa\leq\lambda<\mu$ are infinite cardinals then

\begin{center}
$\lnot\mathsf{KH}(\kappa,\lambda)\land\mathsf{KH}(\kappa,\mu)\rightarrow\mathsf{KH}(\lambda^+,\mu)$
\end{center}

In particular, letting $\kappa=\lambda=\aleph_1$ and $\mu=\aleph_2$ we see that $\lnot\mathsf{KH}(\aleph_1,\aleph_1)\land\mathsf{KH}(\aleph_1,\aleph_2)\land\lnot\mathsf{KH}(\aleph_2,\aleph_2)$ is inconsistent. One of our main goals is to show that all other seven combinations of $\mathsf{KH}(\aleph_1,\aleph_1)$, $\mathsf{KH}(\aleph_1,\aleph_2)$ and $\mathsf{KH}(\aleph_2,\aleph_2)$ are in fact consistent, relative to large cardinals.

Organization of the paper: Section 1 collects some basic facts about forcing and trees that are needed later. Section 2 defines the two-cardinal Kurepa Hypothesis $\mathsf{KH}(\kappa,\lambda)$ and records some simple properties. Section 3 briefly discusses the case of singular $\lambda$, and its relation to other principles; we show that Chang's conjecture $(\aleph_{\omega+1},\aleph_\omega)\twoheadrightarrow(\aleph_{n+1},\aleph_n)$ implies $\lnot\mathsf{KH}(\aleph_{n+1},\aleph_\omega)$, which in turn implies all scales are bad at cofinality $\aleph_{n+1}$. Section 4 carries out the task of showing all combinations $\mathsf{KH}(\aleph_1,\aleph_1)$, $\mathsf{KH}(\aleph_1,\aleph_2)$ and $\mathsf{KH}(\aleph_2,\aleph_2)$ that are not obviously inconsistent are in fact consistent. Section 5 uses similar forcing techniques to solve a problem about one-cardinal trees; we show for singular $\kappa$ the consistency of ``there is a $\kappa^+$-Kurepa tree and there is no $\kappa^{++}$-Aronszajn tree'', extending a result of Cummings. Section 6 lists some questions.

\section{Preliminaries}

\subsection{Notations}

$\mathcal{P}_\kappa(\lambda)$ means $\{x\subseteq\lambda:|x|<\kappa\}$, and $[A]^\kappa$ means $\{x\subseteq A:|x|=\kappa\}$. We use ${}^X Y$ for the set of functions from $X$ to $Y$, and ${}^{<\kappa}A$ for $\bigcup_{\alpha<\kappa}{}^\alpha A$.

We will frequently use $\mathcal{F}\upharpoonright x$ to denote ``restriction of the family $\mathcal{F}$'', whose exact meaning depends on context. When $\mathcal{F}$ is a family of sets (especially sets of ordinals), $\mathcal{F}\upharpoonright x$ denotes $\{A\cap x:A\in\mathcal{F}\}$. When $\mathcal{F}$ is a family of functions, $\mathcal{F}\upharpoonright x$ denotes $\{s\upharpoonright x:s\in\mathcal{F}\}$.

\subsection{Forcing}

We force downward, so $p\leq q$ means $p$ is a strengthening of $q$. A forcing poset is a pre-order $(\mathbb{P},\leq)$ with a distinguished maximal element $1_{\mathbb{P}}$. We often identify the poset with its underlying set $\mathbb{P}$. We say ``$\mathbb{P}$ forces $\varphi$'' if the maximal element $1_{\mathbb{P}}$ forces $\varphi$.

We often use blackboard bold letters such as $\mathbb{P}$, $\mathbb{Q}$, $\mathbb{R}$ to denote forcing posets, and the corresponding roman letters $P$, $Q$, $R$ to denote generic filters. Whether $P$ means an arbitrary generic or a particular one will depend on context. For example, we might say ``since $\mathbb{P}$ is $\kappa^+$-distributive, $V[P]$ has the same subsets of $\kappa$ as $V$''; here $P$ is understood to be an arbitrary generic filter over $V$. On the other hand, if $\mathbb{P}$ projects onto some $\mathbb{P}_\alpha$ and we say ``any branch of $T$ in $V[P]$ is already in $V[P_\alpha]$'', then $P$ means an arbitrary $\mathbb{P}$-generic while $P_\alpha$ means the $\mathbb{P}_\alpha$-generic induced from $P$.

A forcing poset $\mathbb{P}$ is $\kappa$-closed if any decreasing sequence of length strictly less than $\kappa$ has a lower bound; $\mathbb{P}$ is $\kappa$-distributive if it does not add new sequences of ordinals of length strictly less than $\kappa$; $\mathbb{P}$ is $\kappa$-directed closed if any directed subset of size strictly less than $\kappa$ has a lower bound.

\iffalse
\begin{lemma}\label{forcing_closed_preserve_closed}
If $\mathbb{P}$ and $\mathbb{Q}$ are $\kappa$-closed forcings, then $\mathbb{P}$ remains $\kappa$-closed in $V[Q]$ and vice versa.
\end{lemma}

\begin{lemma}[$\Delta$-system lemma]\label{delta_system_lemma}
Suppose $\omega\leq\kappa<\lambda$ are regular cardinals such that $\forall\theta<\lambda,\ \theta^{<\kappa}<\lambda$, and $(A_\alpha:\alpha<\lambda)$ is a sequence of sets such that $|A_\alpha|<\kappa$ for each $\alpha$. Then there exists a set $R$ (called the root) and some $I\in[\lambda]^\lambda$ such that $A_\alpha\cap A_\beta=R$ for any distinct $\alpha,\beta\in I$.
\end{lemma}
\fi

\begin{lemma}[Easton]\label{forcing_Easton}
Assume $\kappa$ is a regular cardinal, $\mathbb{P}$ and $\mathbb{Q}$ are forcing posets, $\mathbb{P}$ is $\kappa$-cc and $\mathbb{Q}$ is $\kappa$-closed. Suppose $P$ is $\mathbb{P}$-generic over $V$ and $Q$ is $\mathbb{Q}$-generic over $V$.

(i) $\mathbb{P}$ is $\kappa$-cc in $V[Q]$. Therefore, any antichain of $\mathbb{P}$ in $V[Q]$ is already in $V$, and $P$ is generic over $V[Q]$.

(ii) $\mathbb{Q}$ is $\kappa$-distributive in $V[P]$, and $Q$ is generic over $V[P]$.
\end{lemma}

A \textit{projection} from $\mathbb{Q}$ onto $\mathbb{P}$ is a map $\pi:\mathbb{Q}\twoheadrightarrow\mathbb{P}$ that preserves order and the distinguished maximal element, such that $\forall q\forall p\leq\pi(q)\exists q'\leq q[\pi(q')\leq p]$. Often times we can actually find $q'$ with $\pi(q')=p$, which implies $\pi$ is surjective, but in general it may not be surjective, although the image is always dense in $\mathbb{P}$. Also, the preimage of a dense open set is dense open, so any $\mathbb{Q}$-generic induces a $\mathbb{P}$-generic. If $P$ is a $\mathbb{P}$-generic, then in $V[P]$ we can consider the \textit{quotient} $\mathbb{Q}/P:=\{q\in\mathbb{Q}:\pi(q)\in P\}$, ordered as a subposet of $\mathbb{Q}$. It can be shown that $\mathbb{Q}*(\mathbb{Q}/\dot{P})$ is forcing equivalent to $\mathbb{Q}$.

$\mathbb{P}$ is said to be \textit{$\kappa$-Knaster} if for any $A\in[\mathbb{P}]^\kappa$, there exists a $B\in [A]^\kappa$ whose elements are pairwise compatible. There are many situations where $\kappa$-Knaster comes more handy than mere $\kappa$-cc.

\begin{lemma}\label{forcing_Knaster}
(i) If $\mathbb{P}$ is $\kappa$-Knaster and $1_{\mathbb{P}}\Vdash\dot{\mathbb{Q}}$ is $\kappa$-Knaster, then $\mathbb{P}*\dot{\mathbb{Q}}$ is $\kappa$-Knaster.

(ii) If $\mathbb{P}$ and $\mathbb{Q}$ are $\kappa$-Knaster then so is $\mathbb{P}\times\mathbb{Q}$. In particular, $\mathbb{P}^n$ is $\kappa$-Knaster for every $n$.

(iii) If $\mathbb{Q}$ is $\kappa$-Knaster, $\pi:\mathbb{Q}\rightarrow\mathbb{P}$ is a projection and $P$ is $\mathbb{P}$-generic over $V$, then in $V[P]$ the quotient $\mathbb{Q}/P$ is square-$\kappa$-cc, namely $(\mathbb{Q}/P)^2$ is $\kappa$-cc.
\end{lemma}
\begin{proof}
For (iii), first note that $\mathbb{P}$ is $\kappa$-Knaster, so $\mathbb{P}\times\mathbb{Q}^2$ is $\kappa$-Knaster, and in particular $\mathbb{Q}^2$ is $\kappa$-cc in $V[P]$; in fact it can be proven to be $\kappa$-Knaster. $(\mathbb{Q}/P)^2$ is a subposet of the $\kappa$-cc poset $\mathbb{Q}^2$, and is thus $\kappa$-cc.
\end{proof}

\subsection{Trees}

A \textit{tree} is a partial order $(T,\preceq)$ such that for any $t\in T$, $\{s\in T:s\prec t\}$ is well-orderd by $\preceq$. An element $t\in T$ is also referred to as a node, and its \textit{rank} is the order type of $\{s\in T:s\prec t\}$. The $\alpha$-th \textit{level} of $T$ is the set of nodes with rank $\alpha$, denoted $\mathcal{L}_\alpha(T)$, and the \textit{height} of $T$ is the least $\alpha$ such that $\mathcal{L}_\alpha(T)=\emptyset$. A \textit{branch} of $T$ is a maximal linear subset; a \textit{cofinal branch} is a branch that has the same height as the $T$; in this paper we call a cofinal branch simply a branch, although the distinction can matter elsewhere.

A $\kappa$-tree is a tree of size and height $\kappa$. It is often harmless to consider only subtrees of ${}^{<\kappa}\kappa$, but sometimes we do need the important fact that a $\kappa$-tree is a binary relation on $\kappa$, and thus coded by a subset of $\kappa$.

A $\kappa$-tree is \textit{slim} if each level has size strictly below $\kappa$. This is mostly intended for successor $\kappa$; if $\kappa$ is inaccessible one usually require the $\alpha$-th level has size at most $|\alpha|$. 

A $\kappa$-\textit{Aronszajn} tree is a slim $\kappa$-tree without cofinal branch. A $\kappa$-\textit{Suslin} tree is a slim $\kappa$-tree without cofinal branch or antichain of size $\kappa$. A $\kappa$-\textit{Kurepa} tree is a slim $\kappa$-tree with at least $\kappa^+$ many cofinal branches.

For historical reason, the tree property $\mathsf{TP}(\kappa)$ denotes ``there is no $\kappa$-Aronszajn tree'', the Suslin Hypothesis $\mathsf{SH}(\kappa)$ denotes ``there is no $\kappa$-Suslin tree'', while the Kurepa Hypothesis $\mathsf{KH}(\kappa)$ denotes ``there exists a $\kappa$-Kurepa tree''.

A \textit{$\kappa$-Kurepa family} is a collection $\mathcal{F}\subseteq\mathcal{P}(\kappa)$ such that $\forall\alpha<\kappa,\ |\mathcal{F}\upharpoonright\alpha|<\kappa$ and $|\mathcal{F}|\geq\kappa^+$. There is a $\kappa$-Kurepa tree if and only if there is a $\kappa$-Kurepa family: for the forward direction, identify the tree with $\kappa$ and consider the branches, which form a Kurepa family; for the backward direction identify $\mathcal{P}(\kappa)$ with ${}^{\kappa}2$ and consider the tree $\{f\upharpoonright\alpha:f\in\mathcal{F},\ \alpha<\kappa\}\subseteq{}^{<\kappa}2$.

When showing non-existence of Aronszajn or Kurepa trees, we often need \textit{branch lemmas} to ensure that certain forcing extensions do not add new branches to existing trees. Below are two very useful branch lemmas due to Unger \cite{unger2012fragility,unger2015fragility}.

\begin{lemma}\label{forcing_branch_Unger_formerly_closed}
Suppose $\kappa\leq\lambda$ are regular, $\mathbb{P}$ is $\kappa$-cc, $\mathbb{Q}$ is $\kappa$-closed, and $\exists\tau<\kappa,\ 2^\tau\geq\lambda$. If $T$ is a slim $\lambda$-tree in $V[P]$, then $T$ does not gain any new branch in $V[P][Q]$.
\end{lemma}

The case most important for us is when $\kappa=\lambda$ is some successor cardinal, so automatically $\exists\tau<\kappa,\ 2^\tau\geq\lambda$. This special case (for $\kappa=\lambda=\omega_1$) also appears as \cite[Lemma 9]{devlin1973kurepa}.

\begin{lemma}\label{forcing_branch_Unger_square-cc}
If $\mathbb{P}$ is square-$\kappa$-cc, and $T$ is a slim $\kappa$-tree, then $T$ does not gain any new branch in $V[P]$.
\end{lemma}

\section{$\mathsf{KH}(\kappa,\lambda)$}

Recall that if $\mathcal{F}$ is a family of sets of ordinals, $\mathcal{F}\upharpoonright x$ denotes $\{A\cap x:A\in \mathcal{F}\}$.

\begin{definition}
Suppose $\kappa\leq\lambda$ are uncountable cardinals.

(i) A \textit{slim $(\kappa,\lambda)$-family} is a collection $\mathcal{F}\subseteq\mathcal{P}(\lambda)$ such that $|\mathcal{F}\upharpoonright x|<\kappa$ for every $x\in\mathcal{P}_\kappa(\lambda)$.

(ii) $\mathsf{KH}(\kappa,\lambda,\mu)$ is the statement that there exists a slim $(\kappa,\lambda)$-family of cardinality $\mu$.

(iii) $\mathsf{KH}(\kappa,\lambda)$ is $\mathsf{KH}(\kappa,\lambda,\lambda^+)$.
\end{definition}

If $\lambda$ is regular, then since every $x\in\mathcal{P}_\lambda(\lambda)$ is contained in some $\alpha<\lambda$, we see that $\mathsf{KH}(\lambda,\lambda)$ is equivalent to $\mathsf{KH}(\lambda)$. We can interpret a $\lambda$-Kurepa family to be a large family that is ``locally small''; then a slim $(\kappa,\lambda)$-family is also locally small, but ``at a different scale''.

We record some simple properties of two-cardinal Kurepa Hypotheses.

\begin{proposition}
(i) If $\mu'\leq\mu$, then $\mathsf{KH}(\kappa,\lambda,\mu)$ implies $\mathsf{KH}(\kappa,\lambda,\mu')$.

(ii) If $\lambda'\leq\lambda$, then $\mathsf{KH}(\kappa,\lambda',\mu)$ implies $\mathsf{KH}(\kappa,\lambda,\mu)$.

(iii) If $\kappa$ is a strong limit cardinal, then $\mathsf{KH}(\kappa,\lambda,2^\lambda)$ holds for any $\lambda\geq\kappa$.

(iv) $\mathsf{KH}(\kappa,\lambda,\lambda)$ always holds.
\end{proposition}
\begin{proof}
(i) and (ii) follow from definition. For (iii), note that any $\mathcal{F}\subseteq\mathcal{P}(\lambda)$ is a slim $(\kappa,\lambda)$-family. For (iv), consider the family $\mathcal{F}$ of all ordinals below $\lambda$, in other words $\mathcal{F}=\lambda$. If $x\in\mathcal{P}_\kappa(\lambda)$ and $\alpha<\lambda$, then $\alpha\cap x$ is either all of $x$ or $\xi\cap x$ where $\xi$ is the least ordinal in $x$ such that $\xi\geq\alpha$. So $|\mathcal{F}\upharpoonright x|\leq|x|+1<\kappa$.
\end{proof}

Because of (iii) above, it is also natural to make the following definition. 

\begin{definition}\label{KH*}
Suppose $\kappa\leq\lambda$ are uncountable cardinals.

(i) A \textit{thin $(\kappa,\lambda)$-family} is a collection $\mathcal{F}\subseteq\mathcal{P}(\lambda)$ such that $|\mathcal{F}\upharpoonright x|\leq|x|$ for every infinite $x\in\mathcal{P}_\kappa(\lambda)$.

(ii) $\mathsf{KH}^*(\kappa,\lambda,\mu)$ is the statement that there exists a thin $(\kappa,\lambda)$-family of cardinality $\mu$.

(iii) $\mathsf{KH}^*(\kappa,\lambda)$ is $\mathsf{KH}^*(\kappa,\lambda,\lambda^+)$.
\end{definition}

Note that $\mathsf{KH}^*(\kappa,\lambda)$ is meaningful for inaccessible as well as singular $\kappa$. The only occassion where $\mathsf{KH}^*(\kappa,\lambda)$ will be relevant for us is Theorem \ref{KH_L}. Another interesting possibility is $\kappa=\omega$, which is best saved for a future paper.

The following simple fact motivated much of the current paper.

\begin{lemma}\label{KH_constraint}
If $\kappa\leq\lambda<\mu$ are infinite cardinals, then $\lnot\mathsf{KH}(\kappa,\lambda)\land
\mathsf{KH}(\kappa,\mu)\rightarrow\mathsf{KH}(\lambda^+,\mu)$.
\end{lemma}
\begin{proof}
By $\mathsf{KH}(\kappa,\mu)$ there exists $\mathcal{F}\subseteq\mathcal{P}(\mu)$ with $|\mathcal{F}|=\mu^+$ and $|\mathcal{F}\upharpoonright x|<\kappa$ for any $x\in\mathcal{P}_\kappa(\mu)$. We claim that $\mathcal{F}$ also witnesses $\mathsf{KH}(\lambda^+,\mu)$, namely $|\mathcal{F}\upharpoonright y|\leq\lambda$ for any $y\in\mathcal{P}_{\lambda^+}(\mu)$. If there were a $y$ for which this fails, without loss of generality $|y|=\lambda$, then since $|(\mathcal{F}\upharpoonright y)\upharpoonright x|=|\mathcal{F}\upharpoonright x|<\kappa$ for each $x\in\mathcal{P}_\kappa(y)$, identifying $y$ with $\lambda$ we see that $\mathcal{F}\upharpoonright y$ witnesses $\mathsf{KH}(\kappa,\lambda)$, contradicting the assumption.
\end{proof}

Note that $\mathsf{KH}(\kappa,\lambda,\mu)$ is upward absolute in mild extensions.

\begin{lemma}\label{KH_upward}
Suppose $V\subseteq W$ are transitive models of $\mathsf{ZFC}$, $\kappa$ is a regular cardinal in $V$, and $\mathcal{P}_\kappa(\lambda)^V$ is cofinal in $\mathcal{P}_\kappa(\lambda)^W$. Then $V\models\mathsf{KH}(\kappa,\lambda,\mu)$ implies $W\models\mathsf{KH}(\kappa,\lambda,\mu)$. In particular, this is true if $W=V[P]$ where $\mathbb{P}$ is either $\kappa$-distributive or $\kappa$-cc.\end{lemma}
\begin{proof}
In $V$ let $\mathcal{F}\subseteq\mathcal{P}(\lambda)$ be a slim $(\kappa,\lambda)$-family with size $\mu$. It suffices to show that $\mathcal{F}$ remains slim in $W$. Given an arbitrary $x\in(\mathcal{P}_\kappa(\lambda))^W$, there exists $y\in(\mathcal{P}_\kappa(\lambda))^V$ such that $x\subseteq y$. Clearly $|\mathcal{F}\upharpoonright x|^W\leq|\mathcal{F}\upharpoonright y|^W\leq|\mathcal{F}\upharpoonright y|^V<\kappa$.

If $\mathbb{P}$ is $\kappa$-distributive, then $\mathcal{P}_\kappa(\lambda)^V=\mathcal{P}_\kappa(\lambda)^W$, hence the upward absoluteness. If $\mathbb{P}$ is $\kappa$-cc, then by a standard covering argument, $\mathcal{P}_\kappa(\lambda)^V$ is cofinal in $\mathcal{P}_\kappa(\lambda)^W$.
\end{proof}

Be aware that under the same assumption, $V\models\mathsf{KH}(\kappa,\lambda)$ does not necessarily imply $W\models\mathsf{KH}(\kappa,\lambda)$, since $(\lambda^+)^V$ may be collapsed in $W$. We do get $W\models\mathsf{KH}(\kappa,\lambda)$ if $V\models\mathsf{KH}(\kappa,\lambda,\mu)$, $\mu>\lambda$ and $\mu$ is preserved in $W$.

Recall that the \textit{Chang's Conjecture} $(\lambda^+,\lambda)\twoheadrightarrow(\kappa^+,\kappa)$ says if $\mathcal{L}$ is any first order language with a distinguished unary predicate $U$, and $\mathcal{M}=(M,U^\mathcal{M},\dots)$ is an $\mathcal{L}$-structure with $|M|=\lambda^+$ and $|U^\mathcal{M}|=\lambda$, then $\mathcal{M}$ has an elementary submodel $\mathcal{N}=(N,U^\mathcal{N},\dots)$ with $|N|=\kappa^+$ and $|U^\mathcal{N}|=\kappa$. The following fact is probably folklore; a version of it appears as \cite[Theorem 3.1]{golshani2019generalized} as well as \cite[Remark 7.6.39]{todorcevic2007walks}.

\begin{lemma}\label{Chang_imply_KH}
$(\lambda^+,\lambda)\twoheadrightarrow(\kappa^+,\kappa)$ implies $\lnot\mathsf{KH}(\kappa^+,\lambda)$.
\end{lemma}
\begin{proof}
Suppose for contradiction that $\mathcal{F}$ is a slim $(\kappa^+,\lambda)$-family and $|\mathcal{F}|=\lambda^+$. Consider the bipartite graph on the disjoint union $\mathcal{F}\sqcup\lambda$ where $A\subseteq\lambda$ is connected to $\alpha<\lambda$ if and only if $\alpha\in A$; interpret the unary predicate $U$ as $\lambda$. An elementary sub-model $\mathcal{G}\sqcup x$ of type $(\kappa^+,\kappa)$ would imply $|\mathcal{F}\upharpoonright x|\geq\kappa^+$, contradicting that $\mathcal{F}$ is slim.
\end{proof}

\section{$\mathsf{KH}(\kappa,\lambda)$ for singular $\lambda$}

Note that $\mathsf{KH}(\kappa,\lambda)$ makes perfect sense for singular $\lambda$. By Lemma \ref{Chang_imply_KH}, we know that if $(\aleph_{\omega+1},\aleph_\omega)\twoheadrightarrow(\aleph_{n+1},\aleph_n)$ holds then $\mathsf{KH}(\aleph_{n+1},\aleph_\omega)$ fails. 

The consistency of $(\aleph_{\omega+1},\aleph_\omega)\twoheadrightarrow(\aleph_{1},\aleph_0)$ (and thus $\lnot\mathsf{KH}(\aleph_1,\aleph_\omega)$) was established by Levinski, Magidor and Shelah \cite{levinski1990chang} using an assumption around the level of huge cardinal; Hayut \cite{hayut2017chang} reduced the large cardinal needed to less than a supercompact cardinal. The consistency of $(\aleph_{\omega+1},\aleph_\omega)\twoheadrightarrow(\aleph_{n+1},\aleph_n)$ for $n=1,2$ is a notorious open problem, while for $n\geq 3$ it is known to be inconsistent. The reason is that $(\aleph_{\omega+1},\aleph_\omega)\twoheadrightarrow(\aleph_{n+1},\aleph_n)$ implies all scales on $\aleph_\omega$ are bad at cofinality $\aleph_{n+1}$, which seems currently out of reach for $n=1,2$, and is outright false for $n\geq 3$ by pcf theory (this was pointed out in \cite{sharon-viale2010}; see \cite{cummings_notes} for a detailed proof).

One might wonder if $\lnot\mathsf{KH}(\aleph_{n+1},\aleph_\omega)$ is nevertheless consistent for $n\geq 3$, or if it is more tractable than the corresponding Chang's Conjecture for $n=1,2$. We show that it is not: it turns out that $\lnot\mathsf{KH}(\aleph_{n+1},\aleph_\omega)$ also implies all scales are bad at cofinality $\aleph_{n+1}$.

Recall that for an infinite $A\subseteq\omega$, a \textit{scale} on $\prod_{k\in A} \aleph_k$ is a sequence $(f_\alpha:\alpha<\aleph_{\omega+1})$ of functions such that:

(i) $f_\alpha\in\prod_{k\in A} \aleph_k$;

(ii) $\alpha<\beta\rightarrow f_\alpha<^* f_\beta$, where $<^*$ means eventual dominance;

(iii) for any $g\in\prod_{k\in A} \aleph_k$, there is $\alpha<\aleph_{\omega+1}$ such that $g<^* f_\alpha$.

As shown by Shelah, there always exists some $A$ that carries a scale.

Given a scale on $\prod_{k\in A} \aleph_k$, an ordinal $\alpha<\aleph_{\omega+1}$ with $\mathrm{cf}(\alpha)>\omega$ is \textit{good} if it satisfies any of the following equivalent conditions.

(a) There exist an unbounded subset $C\subseteq\alpha$ and $k_0<\omega$ such that for all $\beta<\eta$ in $C$ and all $k\in A\setminus k_0$, we have $f_\beta(k)<f_\eta(k)$.

(b) For every unbounded $C\subseteq\alpha$, there exist an unbounded $D\subseteq C$ and $k_0<\omega$ such that for all $\beta<\eta$ in $D$ and all $k\in A\setminus k_0$, we have $f_\beta(k)<f_\eta(k)$.

(c) There exists a sequence $\langle h_i:i<\mathrm{cf}(\alpha)\rangle$ of ordinal functions on $A$ such that: 
	
(c-1) $\langle h_i:i<\mathrm{cf}(\alpha)\rangle$ is pointwise increasing, i.e., if $i<j<\mathrm{cf}(\alpha)$ then $h_i(k)\leq h_j(k)$ for all $k\in A$;
	
(c-2) $\langle f_\beta:\beta<\alpha\rangle$ is cofinally interleaved with $\langle h_i:i<\mathrm{cf}(\alpha)\rangle$ in $<^*$, i.e., $\forall\beta<\alpha\exists i<\mathrm{cf}(\alpha)\ f_\beta<^* h_i$ and $\forall i<\mathrm{cf}(\alpha)\exists\beta<\alpha \ h_i<^* f_\beta$.

A scale is \textit{good} if it is good almost everywhere; more precisely, there is a club $C\subseteq\aleph_{\omega+1}$ such that every point in $C\cap\{\alpha<\aleph_{\omega+1}:\mathrm{cf}(\alpha)>\omega\}$ is good.

Note that if there is a good scale then there is a scale that is good at \textit{every} ordinal $\alpha$ with $\mathrm{cf}(\alpha)>\omega$. Indeed, let $(f_\alpha:\alpha<\aleph_{\omega+1})\subseteq\prod_{k\in A}\aleph_k$ be a good scale and $C\subseteq\aleph_{\omega+1}$ be a club that witnesses the scale is good almost everywhere. Enumerate $C$ as $\{\alpha_\xi:\xi<\aleph_{\omega+1}\}$ and define another scale by $g_\xi:=f_{\alpha_\xi}$. Condition (c) easily implies that in the new scale, every $\alpha<\aleph_{\omega+1}$ with uncountable cofinality is good. By the same argument, if there is a scale that is good at almost all points of cofinality $\aleph_n$, then there is another scale that is good at all points of cofinality $\aleph_n$.

\begin{theorem}
Let $n<\omega$. If there is a scale on $\aleph_\omega$ that is good at all points of cofinality $\aleph_{n+1}$, then $\mathsf{KH}(\aleph_{n+1},\aleph_\omega)$ holds.
\end{theorem}
\begin{proof}
Let $(f_\alpha:\alpha<\aleph_{\omega+1})\subseteq\prod_{k\in A}\aleph_k$ be such a scale. By condition (c), if we truncate the scale at some $\alpha<\aleph_{\omega+1}$ and change each $f_\alpha$ on finitely many values, the resulting scale is still good at all points of cofinality $\aleph_{n+1}$. Thus we may assume that $\aleph_{k-1}\leq f_\alpha(k)<\aleph_k$ for every $\alpha<\aleph_{\omega+1}$ and $k\in A$. In particular, each $f_\alpha$ is a strictly increasing sequence, so we may identify it with a countable subset of $\aleph_\omega$. We claim that they witness $\mathsf{KH}(\aleph_{n+1},\aleph_\omega)$.

Let $X\subseteq\aleph_{\omega}$ be of size $\aleph_{n}$. We want to show that $|\{f_\alpha\cap X:\alpha<\aleph_{\omega+1}\}|\leq\aleph_{n}$. Otherwise there exists $P\subseteq\aleph_{\omega+1}$ of size $\aleph_{n+1}$ such that the sets $f_\alpha\cap X,\ \alpha\in P$ are distinct. Replacing $P$ by a subset of order type $\aleph_{n+1}$, we get a subsequence $(f_{\alpha_i}:i<\aleph_{n+1})$ of the scale such that $f_{\alpha_i}\cap X,\ i<\aleph_{n+1}$ are distinct. Since $\alpha=\sup\alpha_i$ is good, using condition (b) we can pass to a further subsequence so that for some finite $B\subseteq A$ and all $k\in A\setminus B$, the sequence $(f_{\alpha_i}(k):i<\aleph_{n+1})$ is strictly increasing.

For $k\in A$ let $X_k=X\cap[\aleph_{k-1},\aleph_k)$. Since $|X_k|\leq|X|=\aleph_{n}$, if $k\in A\setminus B$ then $(f_{\alpha_i}(k):i<\aleph_{n+1})$ is eventually disjoint from $X_k$. By taking $i$ to be large enough we can make this uniform in $k\in A\setminus B$. For each finitely many $k\in B$, $f_{\alpha_i}\cap X_k$ has at most $|X_k|\leq\aleph_{n}$ many possibilities. Putting these together, we see that $f_{\alpha_i}\cap\bigcup_{k\in A}X_k=f_{\alpha_i}\cap X$ has at most $\aleph_{n}$ many possibilities, a contradiction.
\end{proof}

Taking contrapositive, we see that $\lnot\mathsf{KH}(\aleph_{n+1},\aleph_\omega)$ implies every scale has stationarily many bad points of cofinality $\aleph_{n+1}$. Therefore, $\lnot\mathsf{KH}(\aleph_{n+1},\aleph_\omega)$ is inconsistent for $n\geq 3$, and presumably as difficult as Chang's Conjecture for $n=1,2$.

We remark that $\lnot\mathsf{SCH}_{\aleph_{\omega}}$ implies the existence of a good scale \cite[Theorem 19.3]{cummings_notes}, so $\lnot\mathsf{KH}(\aleph_{n+1},\aleph_\omega)$ implies $\mathsf{SCH}_{\aleph_{\omega}}$.

\section{Combinations of $\mathsf{KH}(\kappa,\lambda)$}

Recall from Lemma \ref{KH_constraint} that $\lnot\mathsf{KH}(\aleph_1,\aleph_1)\land\mathsf{KH}(\aleph_1,\aleph_2)\land\lnot\mathsf{KH}(\aleph_2,\aleph_2)$ is inconsistent. In this section we shall show the consistency of all other seven combinations of $\mathsf{KH}(\aleph_1,\aleph_1)$, $\mathsf{KH}(\aleph_1,\aleph_2)$ and $\mathsf{KH}(\aleph_2,\aleph_2)$. 

We need both tools to destroy Kurepa families and tools to create them. Our main way of destroying Kurepa families is a two-cardinal version of Silver's branch lemma. To this end, we need to reformulate Kurepa families in terms of trees.

\begin{definition}
(i) A \textit{$(\kappa,\lambda)$-tree} is a collection $\mathcal{F}=(\mathcal{F}_x:x\in\mathcal{P}_\kappa(\lambda))$ such that $\mathcal{F}_x\subseteq{}^x2$ is a nonempty set of functions from $x$ to $2=\{0,1\}$, and $\mathcal{F}$ is downward closed in the sense that if $x\subseteq y$ and $s\in\mathcal{F}_y$ then $s\upharpoonright x\in\mathcal{F}_x$.

(ii) $\mathcal{F}$ is \textit{slim} if $\forall x\in\mathcal{P}_\kappa(\lambda),\ |\mathcal{F}_x|<\kappa$.

(iii) A \textit{branch} through $\mathcal{F}$ is a function $b:\lambda\rightarrow 2$ such that $b\upharpoonright x\in\mathcal{F}_x$ for every $x\in\mathcal{P}_\kappa(\lambda)$.
\end{definition}

Clearly $\mathsf{KH}(\kappa,\lambda,\mu)$ holds if and only if there is a slim $(\kappa,\lambda)$-tree with $\mu$ branches. One advantage of the tree formulation is as follows. A slim $(\kappa^+,\lambda)$-trees is coded by an element of $\prod_{x\in\mathcal{P}_{\kappa^+}(\lambda)}{}^\kappa({}^x2)$, in other words a subset of $\bigsqcup_{x\in\mathcal{P}_{\kappa^+}(\lambda)}\kappa\times x$. Thus if $\lambda^\kappa=\lambda$ then a slim $(\kappa^+,\lambda)$-tree is coded by a subset of $\lambda$. Consequently, if $V\subseteq W$ is a forcing extension and $\mathcal{P}_{\kappa^+}(\lambda)^V=\mathcal{P}_{\kappa^+}(\lambda)^W$, then any slim $(\kappa^+,\lambda)$-tree in $W$ already appears in some small intermediate extension, and we are often able to argue the quotient forcing does not add branch. Similar argument works if $\mathcal{P}_{\kappa^+}(\lambda)^V$ is cofinal in $\mathcal{P}_{\kappa^+}(\lambda)^W$.

\begin{lemma}\label{KH_Silver}
(i) If $\mathcal{F}$ is a slim $(\tau^+,\lambda)$-tree and $\mathbb{P}$ is a $\tau^+$-closed forcing poset, then $\mathbb{P}$ does not add new branch to $\mathcal{F}$.

(ii) If $\kappa$ is regular and $\lambda>\kappa$ is inaccessible, then $\mathrm{Col}(\kappa,<\lambda)$ forces $\lnot\mathsf{KH}(\tau^+,\kappa)$ for every infinite cardinal $\tau<\kappa$.
\end{lemma}
\begin{proof}
(i) Otherwise, suppose $\dot{b}:\lambda\rightarrow2$ is forced to be a new branch. Let $\theta\leq\tau$ be the least cardinal with $2^\theta>\tau$, so $2^{<\theta}\leq\tau$. Inductively build a tree $(p_s:s\in {}^{<\theta}2)$ of decreasing conditions (i.e., $s\sqsubseteq t\rightarrow p_s\geq p_t$) and a tree $(\alpha_s:s\in {}^{<\theta}2)$ of ordinals in $\lambda$, so that $p_{s^\smallfrown i}\Vdash\dot{b}(\alpha_s)=i$. For each $f\in {}^\theta2$ let $p_f$ be a lower bound of $\{p_{f\upharpoonright\xi}:\xi<\theta\}$, and let $x=\{\alpha_s:s\in {}^{<\theta}2\}$, so $x\in\mathcal{P}_{\tau^+}(\lambda)$. Then for different $f,g\in {}^\theta2$, $p_f$ and $p_g$ decide $\dot{b}\upharpoonright x$ in different ways, contradicting $|\mathcal{F}_x|\leq\tau$.

(ii) Let $\mathbb{C}=\mathrm{Col}(\kappa,<\lambda)$ and $\mathbb{C}_\alpha=\mathrm{Col}(\kappa,<\alpha)$. Fix an infinite $\tau<\kappa$. Replacing $V$ with $V[C_\alpha]$ for some large enough $\alpha$ we may assume $\kappa^\tau=\kappa$, since $(\kappa^\tau)^V<\lambda$ is collapsed to $\kappa$ in $V[C_\alpha]$, and the quotient $\mathbb{C}/C_\alpha\simeq\mathrm{Col}(\kappa,[\alpha,\lambda))$ is still Levy collapse.

A slim $(\tau^+,\kappa)$-tree can be coded by a subset of $\kappa^\tau=\kappa$. Since $\mathbb{C}$ is $\lambda$-cc, any subset of $\kappa$, and thus any $(\tau^+,\kappa)$-tree $\mathcal{F}\in V[C]$ already belongs to some intermediate extension $V[C_\alpha]$. In $V[C_\alpha]$ the tree has at most $2^\kappa\leq(2^{\alpha})^V$ branches, which will be collapsed by the quotient $\mathrm{Col}(\kappa,[\alpha,\kappa))$ to $\kappa$. By (i) the tail does not introduce new branch, and the result follows.
\end{proof}

Let $\kappa<\lambda$ be infinite regular cardinals. We now define a poset that forces $\mathsf{KH}(\kappa^+,\lambda,\mu)$. This is taken almost verbatim from \cite[Lemma 2.1]{friedman-golshani2012independence}, with some minor changes to make the proof of some properties of the poset slightly easier. Also note that what they denote by $\mathsf{KH}(\kappa,\lambda)$ is $\mathsf{KH}(\lambda^+,\kappa)$ in our notation.

First let us explain the idea. There is a standard way to add a $\kappa^+$-Kurepa tree via a $\kappa^+$-closed forcing $\mathbb{Q}$ (it also works for inaccessible cardinal); the conditions are $(T,g)$ where $T\subseteq {}^{<\kappa^+}2$ is a tree of height $\alpha+1$ for some $\alpha<\kappa^+$, each level of size $\leq\kappa$, and $g$ is a partial function from $\kappa^{++}$ to $\mathcal{L}_\alpha(T)$, the set of top level nodes of $T$. It is straightforward to define a variant $\mathbb{Q}_\mu$ that adds $\mu$ many branches, and if $\eta<\mu$ there is a natural projection from $\mathbb{Q}_\mu$ onto $\mathbb{Q}_\eta$. For later use, we would like $\mathbb{Q}_\mu$ to be not only $\kappa^+$-closed, but also \textit{canonically $\kappa$-closed}, i.e., any decreasing sequence $(p_i:i<\tau)$ where $\tau<\kappa$ has a greatest lower bound. This would ensure the quotient poset $\mathbb{Q}_\mu/Q_\eta$ is $\kappa$-closed in $V[Q_\eta]$. Canonical closure can be arranged by further requiring $T$ to be ``complete'' at levels of cofinality below $\kappa$. Namely, if $\alpha<\mathrm{ht}(T)$, $\mathrm{cf}(\alpha)=\tau<\kappa$, and $s\in {}^{\alpha}2$ is such that $s\upharpoonright\beta\in T$ for any $\beta<\alpha$, then $s\in T$. Note that to keep $T$ slim we need to assume $\kappa^{<\kappa}=\kappa$, so in particular $\kappa$ is regular; we do not know if this forcing has a version that works for singular $\kappa$. The forcing in \cite{friedman-golshani2012independence} is a natural generalization of the above idea to two-cardinal setting, with some complications arising from the fact that the ``levels'' of a $(\kappa^+,\lambda)$-tree are not linearly ordered.

Recall our convention that if $\mathcal{F}$ is a family of functions, then $\mathcal{F}\upharpoonright X$ denotes $\{s\upharpoonright X:s\in\mathcal{F}\}$.

\begin{definition}
Let $\kappa<\lambda$ be infinite regular and $\mu$ be arbitrary. Assume $\kappa^{<\kappa}=\kappa$, $2^\kappa=\kappa^+$ (this will be used to estimate chain conditon) and $\lambda^\kappa=\lambda$. Fix a bijection $J:\mathcal{P}_{\kappa^+}(\lambda)\rightarrow\lambda$, and let

\begin{center}
$\mathcal{P}'_{\kappa^+}(\lambda)=\{X\in\mathcal{P}_{\kappa^+}(\lambda):\forall Y\ J(Y)\in X\rightarrow Y\subseteq X\}$.
\end{center}

Note that this is a club in $\mathcal{P}_{\kappa^+}(\lambda)$. The conditions of $\mathbb{Q}(\kappa^+,\lambda,\mu)$ are triples $(X,\mathcal{F},g)$ where:

(a) $X\in\mathcal{P}'_{\kappa^+}(\lambda)$;

(b) $\mathcal{F}\subseteq {}^{X}2$ is a set of functions from $X$ to $2=\{0,1\}$, and $|\mathcal{F}|\leq\kappa$;

(c) $g:\mu\rightharpoonup\mathcal{F}$ is a partial function, and $|g|\leq\kappa$;

(d) (if $\kappa>\omega$) $\mathcal{F}$ is $\kappa$-closed in the following sense. Suppose $\tau<\kappa$ is regular, $(X_i:i<\tau)$ is a strictly increasing sequence in $\mathcal{P}'_{\kappa^+}(\lambda)$ such that $J(X_i)\in X_{i+1}$ and $\bigcup_{i<\tau}X_i=X$; also suppose $f\in {}^{X}2$ is such that $f\upharpoonright X_i\in\mathcal{F}\upharpoonright X_i$ for every $i<\tau$; then $f\in\mathcal{F}$.

Define $(X_1,\mathcal{F}_1,g_1)\leq(X_2,\mathcal{F}_2,g_2)$ if and only if:

(e) either $X_1=X_2$ or $J(X_2)\in X_1$, which implies $X_1\supseteq X_2$;

(f) $\mathcal{F}_1\upharpoonright X_2=\mathcal{F}_2$;

(g) $\mathrm{dom}(g_1)\supseteq\mathrm{dom}(g_2)$ and $\forall\gamma\in\mathrm{dom}(g_2)$, $g_1(\gamma)\supseteq g_2(\gamma)$, i.e., $g_1(\gamma)$ extends $g_2(\gamma)$ as a function.
\end{definition}

\begin{remark}
\label{KH_Q_remark}
(i) In analogy with the one-cardinal case, we can think of $X$ as the top level, and $\mathcal{F}$ as the top level nodes, which determines the part of our final $(\kappa^+,\lambda)$-tree up to the $X$-th level; note that (f) implies the tree is ``pruned''. When $\kappa^+=\lambda$, this is basically the standard forcing adding a $\kappa^+$-Kurepa tree.

Note that if $(X_1,\mathcal{F}_1,g_1)\leq(X_2,\mathcal{F}_2,g_2)$ and $X_1=X_2$, then (f) implies $\mathcal{F}_1=\mathcal{F}_2$, and if $\gamma\in\mathrm{dom}(g_2)$ then $g_1(\gamma)=g_2(\gamma)$.

(ii) The $\tau$ in (d) is unique if it exists. Suppose $(X_i:i<\tau)$ and $(Y_j:j<\tau')$ are as in (d). Then since $J(X_i)\in X_{i+1}\subseteq X=\bigcup_{j<\tau'}Y_j$, we have $J(X_i)\in Y_j$ and hence $X_i\subseteq Y_j$ for some $j$. Similarly each $Y_j$ is contained in some $X_i$, so the two sequences are interleaved and strictly increasing. Then we must have $\tau=\tau'$ since they are regular.

(iii) For any $\beta\in\lambda$, any condition $(X,\mathcal{F},g)$ can be extended to $(X',\mathcal{F}',g')$ so that $X'\ni\beta$, as follows. Let $X'$ be the union of a continuous and strictly increasing sequence $(X_i:i<\kappa)$ in $\mathcal{P}'_{\kappa^+}(\lambda)$ where $X_0\supseteq X\cup\{\beta\}$. Then we can define $\mathcal{F}'$ and $g'$ arbitrarily subject to (f) and (g). Condition (d) is vacuously true because there is no relevant sequence $(X_i:i<\tau)$, for reason similar to the previous remark.

By the same argument, any condition can be extended so that its first coordinate contains any prescribed $X\in\mathcal{P}_{\kappa^+}(\lambda)$.

(iv) We do not demand injectivity or surjectivity of $g$, but the set of conditions for which $g$ is bijective is dense, although we will not need this.

\end{remark}

For completeness we include the proofs of basic properties of $\mathbb{Q}(\kappa^+,\lambda,\mu)$ from \cite{friedman-golshani2012independence}, with some details added.

\begin{lemma}\label{KH_Q_closure_Knaster}
(i) $\mathbb{Q}(\kappa^+,\lambda,\mu)$ is canonically $\kappa$-closed, i.e., any decreasing sequence $(p_i:i<\tau<\kappa)$ has an infimum.

(ii) $\mathbb{Q}(\kappa^+,\lambda,\mu)$ is $\kappa^+$-closed.

(iii) $\mathbb{Q}(\kappa^+,\lambda,\mu)$ is $\kappa^{++}$-Knaster.
\end{lemma}
\begin{proof}

(i) The case $\kappa=\omega$ is trivial, so assume $\kappa>\omega$. Without loss of generality $\tau$ is regular. Suppose $p_i=(X_i,\mathcal{F}_i,g_i)$. The case when $(X_i:i<\tau)$ stabilizes is also easy, so assume they do not. Let $X=\bigcup_{i<\tau}X_i$, which belongs to $\mathcal{P}'_{\kappa^+}(\lambda)$ since the latter is a club. Define $g$ by $\mathrm{dom}(g)=\bigcup_{i<\tau}\mathrm{dom}(g_i)$ and $g(\gamma)=\bigcup\{g_i(\gamma):i<\tau,\ \gamma\in\mathrm{dom}(g_i)\}$. Finally let $\mathcal{F}$ be the set of all $f\in {}^{X}2$ such that $\forall i<\tau,\ f\upharpoonright X_i\in\mathcal{F}_i$; we have $|\mathcal{F}|\leq\kappa^\tau=\kappa$. Note that automatically $\mathrm{ran}(g)\subseteq\mathcal{F}$.

We need to show that $\mathcal{F}$ satisfies condition (d) in the definition of $\mathbb{Q}(\kappa^+,\lambda,\mu)$. If $(Y_j:j<\tau')$ is a relevant sequence then by (ii) of Remark \ref{KH_Q_remark} we have $\tau'=\tau$, and $(X_i:i<\tau)$ and $(Y_j:j<\tau)$ are interleaving. So if $f\in{}^{X}2$ is such that $f\upharpoonright Y_j\in\mathcal{F}\upharpoonright Y_j$ for each $j$, then also $f\upharpoonright X_i\in\mathcal{F}\upharpoonright X_i$ for each $i$, so $f\in\mathcal{F}$ as required.

That $(X,\mathcal{F},g)\leq(X_i,\mathcal{F}_i,g_i)$ is true by construction, and (d) ensures it is an infimum.

(ii) By (i) it suffices to show that any decreasing sequence $(p_i:i<\kappa)$ has a lower bound. The case when $(X_i:i<\kappa)$ stabilizes is again easy, so assume they do not. We define $X$ and $g$ as in (i), and let $\mathcal{F}=\mathrm{ran}(g)\cup\{\bar{f}:f\in\mathcal{F}_i,\ i<\kappa\}$, where $\bar{f}$ is defined as below; note that $|\mathcal{F}|\leq\kappa$.

Fix $f\in\mathcal{F}_i$. Inductively construct $(f_j:i\leq j<\kappa)$ where $f_i=f$ and $f_{j+1}$ is anything in $\mathcal{F}_{j+1}$ that extends $f_j$. If $k$ is a limit, since $k<\kappa$ and $J(X_j)\in X_k$ for $i\leq j<k$, we know from (d) that the unique $f'$ that extends all the $(f_j:i\leq j<k)$ is in $\mathcal{F}_k$, and we let it be $f_k$. Finally let $\bar{f}=\bigcup_{i\leq j<\kappa}f_j$.

Since $(X_i:i<\kappa)$ does not stabilize, $p=(X,\mathcal{F},g)$ satisfies (d) vacuously, and is a lower bound by construction.

(iii) First we claim that $(X_1,\mathcal{F}_1,g_1)\not\perp(X_2,\mathcal{F}_2,g_2)$ if and only if 

\begin{center}
$\mathcal{F}_1\upharpoonright(X_1\cap X_2)=\mathcal{F}_2\upharpoonright(X_1\cap X_2)$, and
\end{center}

\begin{center}
$g_1(\gamma)\upharpoonright(X_1\cap X_2)=g_2(\gamma)\upharpoonright(X_1\cap X_2)$ for all $\gamma\in\mathrm{dom}(g_1)\cap\mathrm{dom}(g_2)$.
\end{center}

Necessity is clear. To see sufficiency, consider $(X,\mathcal{F},g)$ where $X$ is the union of some strictly increasing sequence in $\mathcal{P}'_{\kappa^+}(\lambda)$ and $X\supseteq X_1\cup X_2\cup\{J(X_1),J(X_2)\}$; again, this ensures the condition satisfies (d) vacuously; define $\mathcal{F}$ and $g$ arbitrarily subject to (f) and (g). Then $(X,\mathcal{F},g)$ is a common lower bound.

Suppose $(p_i:i<\kappa^{++})$ are conditions. From $2^\kappa=\kappa^+$ we have $(\kappa^+)^{<\kappa^+}=\kappa^+$, so applying $\Delta$-system lemma twice we get a subfamily of size $\kappa^{++}$ where the $X_i$'s form a $\Delta$-system with root $X$, and the $\mathrm{dom}(g_i)$'s form a $\Delta$-system with root $Z$. There are $|[{}^{X}2]^\kappa|\leq 2^\kappa=\kappa^+$ many possibilities for $\mathcal{F}_i\upharpoonright X$, so we can use pigeonhole principle to get a size $\kappa^{++}$ subfamily where all the $\mathcal{F}_i\upharpoonright X$ equal some $\mathcal{F}^*$. Now for each $\gamma\in Z$ there are $|\mathcal{F}^*|$ many possibilities for $g_i(\gamma)\upharpoonright X$, and $|{}^{Z}(\mathcal{F}^*)|\leq\kappa^\kappa=\kappa^+$, so using pigeonhole principle again we arrive at a size $\kappa^{++}$ subfamily with $g_i(\gamma)\upharpoonright X$ all the same for every $\gamma\in Z$, in other words a family of pariwise compatible conditions.
\end{proof}

We can now show that $\mathbb{Q}(\kappa^+,\lambda,\mu)$ has the intended effect. Be aware that whenever we talk about $\mathbb{Q}(\kappa^+,\lambda,\mu)$, we implicitly assume that $\kappa^{<\kappa}=\kappa$, etc.

\begin{lemma}\label{KH_Q_effect}
$\mathbb{Q}(\kappa^+,\lambda,\mu)$ forces $\mathsf{KH}(\kappa^+,\lambda,\mu)$. In fact it forces $\mathsf{KH}(\kappa^+,\theta,\mu)$ for any $\kappa^+\leq\theta\leq\lambda$, so in particular $2^{\kappa^+}\geq\mu$.
\end{lemma}
\begin{proof}
If $G$ is a generic filter, define $g^*:\mu\rightarrow{}^{\lambda}2$ by $g^*(\gamma)=\bigcup\{g(\gamma):(X,\mathcal{F},g)\in G\}$. By a simple genericity argument, if $\gamma\neq\eta$ then $g^*(\gamma)\neq g^*(\eta)$, and in fact $g^*(\gamma)\upharpoonright A\neq g^*(\eta)\upharpoonright A$ for any $A\in[\lambda]^{\kappa^+}$.

We argue that $\{g^*(\gamma):\gamma<\mu\}$ is a slim $(\kappa^+,\lambda)$-family. By Lemma \ref{KH_Q_closure_Knaster}, $\mathbb{Q}(\kappa^+,\lambda,\mu)$ is $\kappa^+$-closed so does not change $\mathcal{P}_{\kappa^+}(\lambda)$. In $V[G]$ let $X_0\in\mathcal{P}_{\kappa^+}(\lambda)$ be arbitrary. By (iii) of Remark \ref{KH_Q_remark}, there exists some condition $(X,\mathcal{F},g)\in G$ with $X\supseteq X_0$. By the definition of the forcing order, $\{g^*(\gamma)\upharpoonright X:\gamma<\mu\}$ is exactly $\mathcal{F}$, so has size at most $\kappa$; then $\{g^*(\gamma)\upharpoonright X_0:\gamma<\mu\}$ is also bounded in size by $\kappa$, as desired.
\end{proof}

Since $\mathbb{Q}(\kappa^+,\lambda,\mu)$ is $\kappa^+$-closed and $\kappa^{++}$-cc under our cardinal arithmetic assumption, it preserves all cardinals. Thus if $\mu\geq\lambda^+$ then $\mathbb{Q}(\kappa^+,\lambda,\mu)$ forces $\mathsf{KH}(\kappa^+,\lambda)$.

Let $\mathbb{Q}=\mathbb{Q}(\kappa^+,\lambda,\mu)$, and for $\eta<\mu$ let $\mathbb{Q}_\eta=\mathbb{Q}(\kappa^+,\lambda,\eta)$, defined using the same bijection $J:\mathcal{P}_{\kappa^+}(\lambda)\rightarrow\lambda$.

\begin{lemma}\label{KH_Q_projection}
(i) $\pi:\mathbb{Q}\rightarrow\mathbb{Q}_\eta,\ (X,\mathcal{F},g)\mapsto(X,\mathcal{F},g\upharpoonright\eta)$ is a projection. 

(ii) If $(p_i:i<\tau<\kappa)$ is a decreasing sequence then $\pi(\inf p_i)=\inf\pi(p_i)$.

(iii) $\mathbb{Q}/Q_\eta$ is $\kappa$-closed in $V[Q_\eta]$.
\end{lemma}
\begin{proof}
(i) Suppose $(X_1,\mathcal{F}_1,g_1)\in\mathbb{Q}$, $(X_2,\mathcal{F}_2,g_2)\in\mathbb{Q}_\eta$ and $(X_2,\mathcal{F}_2,g_2)\leq(X_1,\mathcal{F}_1,g_1\upharpoonright\eta)$. Consider $(X_2,\mathcal{F}_2,g^*)$ where $\mathrm{dom}(g^*)=\mathrm{dom}(g_1)\cup\mathrm{dom}(g_2)$ and:

(1) If $\gamma\in\mathrm{dom}(g_1)\setminus\mathrm{dom}(g_2)$, in particular $\gamma\geq\eta$, let $g^*(\gamma)$ be anything in $\mathcal{F}_2$ extending $g_1(\gamma)$.

(2) If $\gamma\in\mathrm{dom}(g_2)$, let $g^*(\gamma)=g_2(\gamma)$.

Then $(X_2,\mathcal{F}_2,g^*)$ is a condition below $(X_1,\mathcal{F}_1,g_1)$ in $\mathbb{Q}$ that gets mapped to $(X_2,\mathcal{F}_2,g^*\upharpoonright\eta)=(X_2,\mathcal{F}_2,g_2)$ in $\mathbb{Q}_\eta$. Thus $\pi$ is a projection.

(ii) Recall from Lemma \ref{KH_Q_closure_Knaster} that $\mathbb{Q}$ is canonically $\kappa$-closed. The proof does not involve the $g$ coordinate in a significant way, and it is easily seen that $\pi$ commutes with the infimum.

(iii) In general, if a projection $\pi:\mathbb{Q}\rightarrow\mathbb{P}$ commutes with infima and both posets are canonically $\kappa$-closed, then the quotient is $\kappa$-closed. Indeed, if $P$ is $\mathbb{P}$-generic over $V$, then the quotient is $\mathbb{Q}/P=\{q\in\mathbb{Q}:\pi(q)\in P\}$ ordered as a subposet of $\mathbb{Q}$. If $(q_i:i<\tau<\kappa)$ is a decreasing sequence in $\mathbb{Q}/P$, then $(\pi(q_i):i<\tau<\kappa)$ is a decreasing sequence in $P$; also this sequence belongs to $V$ by closure of $\mathbb{P}$. We want to show $\pi(\inf q_i)=\inf\pi(q_i)$ is in $P$, which would imply $\inf q_i\in\mathbb{Q}/P$ as desired. It suffices to show 

\begin{center}
$\{p\in\mathbb{P}:\forall i<\tau[p\leq \pi(q_i)]\lor\exists i<\tau[p\perp \pi(p_i)]\}$
\end{center}

is a dense set in $V$; then since $P$ cannot contain an element incompatible with any $\pi(q_i)$, it must contain a lower bound of $(\pi(q_i):i<\tau)$, and thus it contains the infimum. The denseness of this set follows either from $\mathbb{P}$ being separated (which is true in our case) or from $\mathbb{P}$ being $\kappa$-closed: given arbitrary $p$, if $p\perp \pi(q_0)$ then we are done; otherwise strengthen $p$ to $p_1$ and compare it with $\pi(q_1)$, etc.
\end{proof}

Another way to get instances of two-cardinal Kurepa Hypotheses is to work over $L$.

\begin{theorem}
\label{KH_L}
$L$ satisfies $\mathsf{KH}(\kappa,\lambda)$ for uncountable regular cardinals $\kappa<\lambda$.
\end{theorem}
\begin{proof}
See \cite[Chapter VII Section 3]{devlin2017constructibility}. Note that what is denoted $\mathsf{KH}(\kappa,\lambda)$ there is $\mathsf{KH}^*(\lambda,\kappa)$ in our notation; cf. Definition \ref{KH*}. It is proved there that $L$ satisfies $\mathsf{KH}^*(\kappa,\lambda)$, which of course implies $\mathsf{KH}(\kappa,\lambda)$.
\end{proof}

For brevity let us denote, e.g., $\mathsf{KH}(\aleph_1,\aleph_1)\land\lnot\mathsf{KH}(\aleph_1,\aleph_2)\land\mathsf{KH}(\aleph_2,\aleph_2)$ by $\mathsf{KH}(\aleph_1,\aleph_1)-\mathsf{KH}(\aleph_1,\aleph_2)+\mathsf{KH}(\aleph_2,\aleph_2)$. Now we have all the tools needed to show that all combinations of $\mathsf{KH}(\aleph_1,\aleph_1)$, $\mathsf{KH}(\aleph_1,\aleph_2)$ and $\mathsf{KH}(\aleph_2,\aleph_2)$ are possible, except for $-\mathsf{KH}(\aleph_1,\aleph_1)+\mathsf{KH}(\aleph_1,\aleph_2)-\mathsf{KH}(\aleph_2,\aleph_2)$, which is ruled out by Lemma \ref{KH_constraint}.

\begin{theorem}\label{KH_comb}
(i) $\mathsf{KH}(\aleph_1,\aleph_1)+\mathsf{KH}(\aleph_1,\aleph_2)+\mathsf{KH}(\aleph_2,\aleph_2)$ is consistent.

(ii) $-\mathsf{KH}(\aleph_1,\aleph_1)+\mathsf{KH}(\aleph_1,\aleph_2)+\mathsf{KH}(\aleph_2,\aleph_2)$ is consistent.

(iii) $\pm\mathsf{KH}(\aleph_1,\aleph_1)-\mathsf{KH}(\aleph_1,\aleph_2)-\mathsf{KH}(\aleph_2,\aleph_2)$ is consistent.

(iv) $\pm\mathsf{KH}(\aleph_1,\aleph_1)-\mathsf{KH}(\aleph_1,\aleph_2)+\mathsf{KH}(\aleph_2,\aleph_2)$ is consistent.

(v) $\mathsf{KH}(\aleph_1,\aleph_1)+\mathsf{KH}(\aleph_1,\aleph_2)-\mathsf{KH}(\aleph_2,\aleph_2)$ is consistent.
\end{theorem}
\begin{proof}
(i) This is true in $L$ by Theorem \ref{KH_L}.

(ii) In $L$ let $\kappa$ be inaccessible, so by Theorem \ref{KH_L} we have $\mathsf{KH}(\aleph_1,\kappa)$; also $\mathsf{KH}(\kappa,\kappa)$ holds trivially since $\kappa$ is a strong limit. Force over $L$ with $\mathrm{Col}(\aleph_1,<\kappa)$, so $\kappa$ becomes $\aleph_2$. By Lemma \ref{KH_upward}, $\mathsf{KH}(\aleph_1,\kappa)$ is preserved since the forcing is $\aleph_1$-closed, and $\mathsf{KH}(\kappa,\kappa)$ is preserved since the forcing is $\kappa$-cc.

(iii) We start with a ground model where the truth value of $\mathsf{KH}(\aleph_1,\aleph_1)$ is already as desired. Let $\kappa$ be inaccessible and force with $\mathrm{Col}(\aleph_2,<\kappa)$. Then $\mathsf{KH}(\aleph_1,\aleph_2)$ and $\mathsf{KH}(\aleph_2,\aleph_2)$ fail by Lemma \ref{KH_Silver}. 

Since the forcing is $\aleph_2$-closed, the status of $\mathsf{KH}(\aleph_1,\aleph_1)$ is as in the ground model $V$. In more detail, if $\mathsf{KH}(\aleph_1,\aleph_1,\aleph_2)$ holds in $V$ then it remains true in the extension by Lemma \ref{KH_upward} since $\aleph_2$ is preserved. If $\mathsf{KH}(\aleph_1,\aleph_1,\aleph_2)$ fails in $V$, in other words $V$ does not have an $\aleph_1$-Kurepa tree, then since the forcing is $\aleph_2$-closed, it adds neither new $\aleph_1$-trees nor new branches to existing $\aleph_1$-trees, so there is still no $\aleph_1$-Kurepa tree in the extension.

(iv) Assume $2^{\aleph_0}=\aleph_1$, $2^{\aleph_1}=\aleph_2$ and $\kappa$ is inaccessible. Let $\mathbb{C}=\mathrm{Col}(\aleph_2,<\kappa)$, $\mathbb{Q}=\mathbb{Q}(\aleph_2,\aleph_2,\kappa)$ and force with $\mathbb{C}\times\mathbb{Q}$. For any $\alpha<\aleph_2$, there are respectively projections $\mathbb{Q}\rightarrow\mathbb{Q}_\alpha:=\mathbb{Q}(\aleph_2,\aleph_2,\alpha)$ and $\mathbb{C}\rightarrow\mathbb{C}_\alpha:=\mathrm{Col}(\aleph_2,<\alpha)$.

From the proof of Lemma \ref{KH_Q_closure_Knaster} we see that $\mathbb{Q}(\kappa^+,\lambda,\mu)$ is $\theta$-Knaster for any $\theta\geq\kappa^{++}$, so here $\mathbb{Q}$ is $\kappa$-Knaster. The Levy collapse $\mathbb{C}$ is also $\kappa$-Knaster, and hence their product. In particular, $\kappa=(\aleph_3)^{V[C\times Q]}$.

Since the forcing is $\aleph_2$-closed, status of $\mathsf{KH}(\aleph_1,\aleph_1)$ is as in the ground model. 

We argue that $\mathsf{KH}(\aleph_1,\aleph_2)$ fails. Any slim $(\aleph_1,\aleph_2)$-tree $\mathcal{F}$ is coded by a subset of $\aleph_2$, so by the chain condition it already appears in some intermediate model $V[C_\alpha\times Q_\alpha]$; in this model $\mathcal{F}$ has at most $2^\alpha<\kappa$ many branches, and $2^\alpha$ is collapsed to $\aleph_2$ by the quotient forcing, so it suffices to show the quotient does not add new branch to $\mathcal{F}$. The quotient forcing is $\mathbb{C}\times\mathbb{Q}/C_\alpha\times Q_\alpha=(\mathbb{C}/C_\alpha)\times(\mathbb{Q}/Q_\alpha)$. $\mathbb{C}/C_\alpha$ is $\aleph_1$-closed (in fact $\aleph_2$-closed) in $V[C_\alpha]$, and $\mathbb{Q}/Q_\alpha$ is $\aleph_1$-closed in $V[Q_\alpha]$ by Lemma \ref{KH_Q_projection}. Since an $\aleph_1$-closed forcing preserves the $\aleph_1$-closure of other posets, $\mathbb{Q}$ remains $\aleph_1$-closed in $V[C_\alpha]$, and $\mathbb{C}/C_\alpha$ remains $\aleph_1$-closed in $V[C_\alpha][Q_\alpha]=V[C_\alpha\times Q_\alpha]$. Similarly $\mathbb{Q}/Q_\alpha$ is $\aleph_1$-closed in $V[C_\alpha\times Q_\alpha]$. We finally conclude that the quotient forcing does not add new branch to $\mathcal{F}$ by Lemma \ref{KH_Silver}.

$\mathsf{KH}(\aleph_2,\aleph_2,\kappa)$ holds in $V[Q]$, and thus in $V[Q][C]$ by Lemma \ref{KH_upward}. Since $\kappa$ is preserved we have $\mathsf{KH}(\aleph_2,\aleph_2)$.

(v) Assume $2^{\aleph_0}=\aleph_1$ and $\kappa$ is inaccessible. Let $\mathbb{C}=\mathrm{Col}(\aleph_2,<\kappa)$, $\mathbb{Q}=\mathbb{Q}(\aleph_1,\aleph_2,\kappa)$ and force with $\mathbb{C}\times\mathbb{Q}$. 

$\mathsf{KH}(\aleph_1,\aleph_1)$ holds because by Lemma \ref{KH_Q_effect}, $\mathbb{Q}(\aleph_1,\aleph_2,\kappa)$ also adds an $\aleph_1$-Kurepa tree. Alternatively it can be seen from Lemma \ref{KH_constraint}.

$\mathsf{KH}(\aleph_1,\aleph_2,\kappa)$ holds in $V[Q]$ and thus in $V[Q][C]$ by Lemma \ref{KH_upward}, because $\mathbb{C}$ is $\aleph_1$-closed in $V[Q]$.

We argue that $\mathsf{KH}(\aleph_2,\aleph_2)$ fails. Any slim $\aleph_2$-tree is coded by a subset of $\aleph_2$, and hence already appears in some intermediate model $V[C_\alpha][Q_\alpha]$ (note that this is not literally true for slim $(\aleph_2,\aleph_2)$-tree), where it has at most $2^\alpha<\kappa$ many branches, and $2^\alpha$ is collapsed to $\aleph_2$ by the quotient forcing. We have the following sequence of extensions:

\begin{center}
$V[C_\alpha][Q_\alpha]\subseteq V[C_\alpha][Q]\subseteq V[C][Q]$
\end{center}

We first argue that the second extension does not add new branch. In $V[C_\alpha]$ the forcing $\mathbb{C}/C_\alpha$ is $\aleph_2$-closed and $\mathbb{Q}$ is $\aleph_2$-cc; the latter follows either from Lemma \ref{forcing_Easton} or from the fact that the definition of $\mathbb{Q}(\kappa^+,\lambda,\mu)$ is unchanged by a $\kappa^+$-closed forcing. Thus we see from Lemma \ref{forcing_branch_Unger_formerly_closed} that no new branch is added.

Finally we need to argue that forcing with $\mathbb{Q}/Q_\alpha$ over $V[C_\alpha][Q_\alpha]$ does not add branch. The definition of $\mathbb{Q}=\mathbb{Q}(\aleph_1,\aleph_2,\kappa)$ is unchanged by $\aleph_1$-closed forcing, so $\mathbb{Q}$ is still $\mathbb{Q}(\aleph_1,\aleph_2,\kappa)$ in $V[C_\alpha][Q_\alpha]$, and in particular still $\aleph_2$-Knaster by Lemma \ref{KH_Q_closure_Knaster}. Then the subposet $\mathbb{Q}/Q_\alpha$ is at least square-$\aleph_2$-cc because $\mathbb{Q}\times\mathbb{Q}$ is $\aleph_2$-cc, so it does not add branch to a slim $\aleph_2$-tree by Lemma \ref{forcing_branch_Unger_square-cc}.
\end{proof}

We note that our use of large cardinals is optimal, since by \cite[Lemma 3.1]{friedman-golshani2012independence}, if $\kappa<\lambda$ are infinite cardinal, $\lambda^\kappa=\kappa$ and $\mathsf{KH}(\kappa^+,\lambda)$ fails, then $\lambda^+$ is inaccessble in $L$.

\section{Aronszajn and Kurepa trees}

In this section we switch gear and come back to one-cardinal trees. Cummings \cite{cummings2018aronszajn} showed that it is consistent to have a regular cardinal $\kappa$ with $\mathsf{TP(\kappa^{++})\land\mathsf{KH}(\kappa^+)}$, namely there is no $\kappa^{++}$-Aronszajn tree but there exists a $\kappa^+$-Kurepa tree. The motivation comes from the fact that there is no $\kappa^+$-Kurepa tree in the standard model of $\mathsf{TP}(\kappa^{++})$, namely the model obtained via the Mitchell forcing $\mathbb{M}(\kappa,\lambda)$, by essentially the same argument as in Lemma \ref{KH_Silver}. Cummings' argument shows that taking the product $\mathbb{M}(\kappa,\lambda)\times\mathbb{Q}(\kappa^+,\kappa^+,\lambda)$ produces a model of $\mathsf{TP(\kappa^{++})\land\mathsf{KH}(\kappa^+)}$; of course most of the work lies in showing there is no $\kappa^{++}$-Aronszajn tree.

We shall show the consistency of $\mathsf{TP(\kappa^{++})\land\mathsf{KH}(\kappa^+)}$ for singular $\kappa$ by combining Cummings' argument with Prikry forcing. In fact, the consistency of $\mathsf{TP(\kappa^{++})\land\lnot\mathsf{KH}(\kappa^+)}$ does not seem immediate either, so we start by demonstrating this.

First we review some standard facts about Mitchell forcing. Suppose $\kappa$ is regular and $\lambda>\kappa$ is at least inaccessible. Denote $\mathbb{P}=\mathrm{Add}(\kappa,\lambda)$ and $\mathbb{P}_\alpha=\mathrm{Add}(\kappa,\alpha)$ for $\alpha<\lambda$. If $p\in\mathbb{P}$, then its natural projection onto $\mathbb{P}_\alpha$ is $p\upharpoonright\alpha$. Assume $\kappa^{<\kappa}=\kappa$, so $\mathbb{P}$ is $\kappa^+$-cc. The conditions of $\mathbb{M}(\kappa,\lambda)$ are $(p,f)$ where 

(i) $p\in\mathrm{Add}(\kappa,\lambda)$;

(ii) $f$ is a partial function on $\lambda$ and $|f|\leq\kappa$;

(iii) for each $\alpha\in\mathrm{dom}(f)$, $f(\alpha)$ is a $\mathbb{P}_\alpha$-name that is forced to be in the poset $\mathrm{Add}(\kappa^+,1)$ as defined in $V[P_\alpha]$, or in symbol $1_{\mathbb{P}_\alpha}\Vdash f(\alpha)\in\dot{\mathrm{Add}}(\kappa^+,1)$.

The order is defined by $(p,f)\leq(q,g)$ if and only if 

(a) $p\leq q$ in $\mathbb{P}$;

(b) $\mathrm{dom}(f)\supseteq\mathrm{dom}(g)$;

(c) for each $\alpha\in\mathrm{dom}(g)$, $p\upharpoonright\alpha\Vdash f(\alpha)\leq g(\alpha)$.

Taken literally, the $\mathbb{M}(\kappa,\lambda)$ defined this way is a proper class. There are multiple ways to circumvent this, such as noticing that for each $\alpha<\lambda$, there are only set many equivalence classes of names $\sigma$ such that $1_{\mathbb{P}_\alpha}\Vdash \sigma\in\dot{\mathrm{Add}}(\kappa^+,1)$.

Let $\mathbb{M}=\mathbb{M}(\kappa,\lambda)$. If $\alpha<\lambda$ is inaccessible, then there is a natural projection from $\mathbb{M}$ onto $\mathbb{M}_\alpha:=\mathbb{M}(\kappa,\alpha)$. We collect some well-known facts about Mitchell forcing.

\begin{fact}
(i) $\mathbb{M}$ is $\kappa$-directed closed and $\lambda$-Knaster.

(ii) There is a $\kappa^+$-directed closed $\mathbb{R}$, called the term space forcing, such that $\mathbb{M}$ is a projection of $\mathbb{P}\times\mathbb{R}$.

(iii) Similarly, if $\alpha<\lambda$ is inaccessible and $M_\alpha$ is $\mathbb{M}_\alpha$-generic over $V$, then there is in $V[M_\alpha]$ an projection $\mathbb{P}^\alpha\times\mathbb{R}^\alpha\twoheadrightarrow\mathbb{M}/M_\alpha$, where $\mathbb{P}^\alpha=\mathrm{Add}(\kappa,[\alpha,\lambda))$ is a Cohen forcing and $\mathbb{R}^\alpha$ is $\kappa^+$-directed closed.
\end{fact}

\begin{theorem}
$\mathsf{TP(\kappa^{++})\land\lnot\mathsf{KH}(\kappa^+)}$ for singular $\kappa$ is consistent relative to large cardinals.
\end{theorem}
\begin{proof}
Let $\kappa$ be indestructibly supercompact and $\lambda>\kappa$ be weakly compact. Let $\mathbb{M}=\mathbb{M}(\kappa,\lambda)$ and $\mathbb{P}\times\mathbb{R}\rightarrow\mathbb{M}$ be the usual projection. We know that $\mathsf{TP}(\kappa^{++})$ holds in $V[M]$; even better, \cite[Theorem 3.2]{honzik2020indestructibility} shows that it is indestructible under any $\kappa^+$-cc forcing that belongs to the Cohen extension; more precisely, for any $\kappa^+$-cc forcing $\overline{\mathbb{P}}\in V[P]$ (by Easton's Lemma \ref{forcing_Easton} it remains $\kappa^+$-cc in $V[P\times R]$ and thus in $V[M]$), if $\overline{P}$ is generic over $V[M]$ then $\mathsf{TP}(\kappa^{++})$ holds in $V[M][\overline{P}]$.

Since $\kappa$ is indestructibly supercompact and $\mathbb{P}$ is $\kappa$-directed closed, $\kappa$ remains supercompact, and in particular measurable, in $V[P]$. There is thus a normal measure $U$ on $\kappa$ in $V[P]$. Since by Easton's Lemma \ref{forcing_Easton} $V[M]$ has the same subsets of $\kappa$ as $V[P]$, we see that $U$ is still a normal measure in $V[M]$. Let $\overline{\mathbb{P}}$ be the Prikry forcing defined with respect to $U$ in $V[P]$, or equivalently in $V[M]$. Since Prikry forcing is $\kappa^+$-cc, by the aforementioned result we have $V[M][\overline{P}]\models\mathsf{TP}(\kappa^{++})$. To finish the proof, we just need to show the following

Claim: $V[M][\overline{P}]\models\mathsf{\lnot KH}(\kappa^+)$.

By \cite[Lemma 5.1]{cummings1998tree}, there are unboundedly many inaccessible $\alpha<\lambda$ for which $U_\alpha:=U\cap\mathcal{P}(\kappa)^{V[P_\alpha]}$ belongs to $V[P_\alpha]$, and thus is a normal measure in $V[P_\alpha]$, or equivalently in $V[M_\alpha]$. Let $\overline{\mathbb{P}}_\alpha$ be the Prikry forcing defined in $V[P_\alpha]$, equivalently in $V[M_\alpha]$, with respect to $U_\alpha$. By Mathias' criterion, if $S$ is the Prikry sequence corresponding to $\overline{P}$, then $S$ induces a filter $\overline{P}_\alpha$ that is $\overline{\mathbb{P}}_\alpha$-generic over $V[M_\alpha]$.

We claim that any slim $\kappa^+$-tree $T$ in $V[M][\overline{P}]$ already belongs to $V[M_\alpha][\overline{P}_\alpha]$ for some inaccessible $\alpha<\lambda$. A condition in $\mathbb{M}$ looks like $(p,f)$, and a condition in $\mathbb{M}*\overline{\mathbb{P}}$ looks like $(p,f,s,\dot{A})$, where $\dot{A}$ is an $\mathbb{M}$-name for a set in $\dot{U}$. Since the forcing is $\lambda$-cc, and $T$ is coded by a subset of $\kappa^+$, its name $\dot{T}$ involves fewer than $\lambda$ many conditions in $\mathbb{M}*\overline{\mathbb{P}}$. By taking $\alpha$ large enough we can ensure all the $(p,f)$'s appearing are in $\mathbb{M}_\alpha$ and all the $\dot{A}$'s are $\mathbb{M}_\alpha$-names, so $\dot{T}$ can be viewed as a $\mathbb{M}_\alpha*\overline{\mathbb{P}}_\alpha$-name.

In $V[M_\alpha]$ we have a projection $\mathbb{P}^\alpha\times\mathbb{R}^\alpha\twoheadrightarrow\mathbb{M}/M_\alpha$ onto the Mitchell remainder. We have the following sequence of extensions.

\begin{center}
$V[M_\alpha][\overline{P}_\alpha]\subseteq V[M_\alpha][P^\alpha][\overline{P}]\subseteq V[M_\alpha][P^\alpha][\overline{P}][R^\alpha]$
\end{center}

Let us explain why $V[M_\alpha][P^\alpha][\overline{P}][R^\alpha]$ makes sense. We assumed $\overline{P}$ to be generic over $V[M]$; by Mathias' criterion it is thus generic over $V[M_\alpha][P^\alpha][R^\alpha]$, which has the same subsets of $\kappa$. But $\overline{\mathbb{P}}$ lives in $V[P]\subseteq V[M_\alpha][P^\alpha]$, so $\overline{P}$ and $R^\alpha$ are mutually generic. Hence it is harmless to switch the order of $\overline{P}$ and $R^\alpha$.

We want to show that neither of these two extensions add new branches to $T$. The first extension is a little subtle to analyze directly; it is an extension via the quotient poset $(\mathbb{P}^\alpha*\overline{\mathbb{P}})/\overline{P}_\alpha$, which is not exactly $\mathbb{P}^\alpha$. But we can argue indirectly as follows. By Lemma \ref{forcing_Knaster} (iii), the quotient poset is square-$\kappa^+$-cc, so does not add branch by Lemma \ref{forcing_branch_Unger_square-cc}.

For the second extension, we apply Lemma \ref{forcing_branch_Unger_formerly_closed} in $V[M_\alpha]$ to the $\kappa^+$-cc forcing $\mathbb{P}^\alpha*\mathbb{\overline{P}}$ and the $\kappa^+$-closed forcing $\mathbb{R}^\alpha$ to see that no new branch is added.

We have shown that $T$ does not gain new branch going from $V[M_\alpha][\overline{P}_\alpha]$ to $V[M][\overline{P}]$. In the former model it has no more than $\alpha$ branches, and $\alpha$ is collapsed to $\kappa^+$ by the Mitchell remainder, so $T$ has at most $\kappa^+$ branches in $V[M][\overline{P}]$.
\end{proof}

\begin{theorem}
$\mathsf{TP(\kappa^{++})\land\mathsf{KH}(\kappa^+)}$ for singular $\kappa$ is consistent relative to large cardinals.
\end{theorem}
\begin{proof}
Let $\kappa$ be indestructibly supercompact, $\lambda>\kappa$ be measurable (this is for simplicity; weakly compact should suffice), $\mathbb{M}=\mathbb{M}(\kappa,\lambda)$, $\mathbb{P}\times\mathbb{R}\twoheadrightarrow\mathbb{M}$ be the usual projection, and $\mathbb{Q}=\mathbb{Q}(\kappa^+,\kappa^+,\lambda)$. Following Cummings, we force with $\mathbb{M}\times\mathbb{Q}$, and then use Prikry forcing to singularize $\kappa$. 

\begin{claim}
Working in $V[M\times Q]$, if $\overline{\mathbb{P}}$ is any $\kappa^+$-cc forcing in the Cohen extension $V[P]$, then $V[M\times Q][\overline{P}]\models\mathsf{TP(\kappa^{++})\land\mathsf{KH}(\kappa^+)}$.
\end{claim}

This in particular applies to the Prikry forcing $\mathbb{\overline{P}}$ as defined in $V[P]$ with respect to some normal measure $U$, and thus would finish the proof. Note that $U$ remains a normal measure in $V[P\times R\times Q]$ and thus in $V[M\times Q]$, because $\mathbb{R}$ and $\mathbb{Q}$ are $\kappa^+$-closed in $V$, so by Easton's Lemma \ref{forcing_Easton} the extension $V[P]\subseteq V[P\times R\times Q]$ does not add $<\kappa^+$-sequences.

We first deal with $\mathsf{KH}(\kappa^+)$, which is the easier part. By Lemma \ref{KH_Q_effect} we have $\mathsf{KH}(\kappa^+,\kappa^+,\lambda)$ in $V[Q]$. We have the following sequence of extensions:

\begin{center}
$V[Q]\subseteq V[R\times Q]\subseteq V[P\times R\times Q]$
\end{center}

The first extension is $\kappa^+$-closed, and the second extension is $\kappa^+$-cc by Easton's Lemma \ref{forcing_Easton}, so by Lemma \ref{KH_upward} we have $\mathsf{KH}(\kappa^+,\kappa^+,\lambda)$ in $V[P\times R\times Q]$, and hence in $V[M\times Q]$ since it contains the witnessing family. Finally $\overline{\mathbb{P}}$ is $\kappa^+$-cc in $V[M\times Q]$, again by Easton, so $\mathsf{KH}(\kappa^+,\kappa^+,\lambda)$ holds in $V[M\times Q][\overline{P}]$.

Next we argue that there is no $\lambda$-Aronszajn tree. Start with an elementary embedding $j:V\rightarrow W$ with critical point $\lambda$. Let $\mathbb{P}^\lambda\times\mathbb{R}^\lambda\twoheadrightarrow j(\mathbb{M})/M$ be the projection in $W[M]$. There is also a projection $j(\mathbb{Q})\twoheadrightarrow\mathbb{Q}$ in $V$ since $\mathbb{Q}=j(\mathbb{Q})_\lambda$. We lift the embedding $j$ successively as follows.

First we force with $\mathbb{P}^\lambda\times \mathbb{R}^\lambda\times j(\mathbb{Q})/Q$ over $V[M\times Q]$. In the resulting model $V[M][P^\lambda\times R^\lambda][j(Q)]$, we can lift the embedding $j:V\rightarrow W$ to

\begin{center}
$j:V[M\times Q]\rightarrow W[j(M)\times j(Q)]$
\end{center}

\noindent where $j(M)$ is the $j(\mathbb{M})$-generic filter induced by $M*(P^\lambda\times R^\lambda)$; it can be checked that $j(M)$ and $j(Q)$ are mutually generic, so the model $W[j(M)\times j(Q)]$ makes sense.

\iffalse
$V[M][Q]=V[Q][M]\subseteq V[Q][M][j(Q)/Q]\subseteq V[Q][M][j(Q)/Q][P^\lambda\times R^\lambda]= V[Q][j(Q)/Q][M][P^\lambda\times R^\lambda]\supseteq V[j(Q)][j(M)]$
\fi

We can now make sense of $j(\mathbb{\overline{P}})$. Note that $j$ restricts to a map from $\overline{\mathbb{P}}$ into $j(\overline{\mathbb{P}})$ that is almost a complete embedding. More precisely, since $V[M\times Q]\models\mathbb{\overline{P}}$ is $\kappa^+$-cc and $\mathrm{crit}(j)=\lambda>\kappa$, if $A\subseteq\overline{\mathbb{P}}$ is a maximal antichain in $V[M\times Q]$ then $j(A)=j[A]$, so $j[A]$ is also maximal.

It follows that if we let $\overline{P}^*$ be $j(\mathbb{\overline{P}})$-generic over $V[M][P^\lambda\times R^\lambda][j(Q)]$, then its preimage under $j$, denoted $\overline{P}$, is $\overline{\mathbb{P}}$-generic over $V[M\times Q]$. Thus we can lift the embedding to

\begin{center}
$j:V[M\times Q][\overline{P}]\rightarrow W[j(M)\times j(Q)][j(\overline{P})]$
\end{center}

\noindent in the model $V[M][P^\lambda\times R^\lambda][j(Q)][j(\overline{P})]$, where $j(\overline{P})=\overline{P}^*$. Let $T$ be a $\lambda$-tree in $V[M\times Q][\overline{P}]$; by standard arguments we may assume it is in $W[M\times Q][\overline{P}]$. The elementary embedding provides a branch of $T$ in $W[j(M)\times j(Q)][j(\overline{P})]$, and we want to show that the branch is already in $W[M\times Q][\overline{P}]$. We have the following sequence of extensions:

\begin{align}
W[M][Q][\overline{P}]&\subseteq W[M][j(Q)][\overline{P}]\\
&\subseteq W[M][P^\lambda][j(Q)][\overline{P}]\\
&\subseteq W[M][P^\lambda][j(Q)][j(\overline{P})]\\
&\subseteq W[M][P^\lambda\times R^\lambda][j(Q)][j(\overline{P})]
\end{align}

It suffices to argue that none of the extension adds branch to $T$.

\iffalse$\overline{P}$ is generic over $V[M\times Q]$, so in particular over $W[M]$. By successive use of Easton, $\overline{\mathbb{P}}$ remains $\kappa^+$-cc in $W[M][j(Q)]$ and $j(\mathbb{Q})$ is $\kappa^+$-distributive in $W[M]$, so it does not add antichains in $\overline{\mathbb{P}}$. Thus $\overline{P}$ is generic over $W[M][j(Q)]$, and we are free to interchange the order of $\overline{P}$ and $j(Q)$.\fi

(1) $j(\mathbb{Q})/Q$ is square-$\lambda$-cc in $W[Q]$ and thus in $W[M][Q][\overline{P}]$, because $\mathbb{M}$ and $\overline{\mathbb{P}}$ are both $\lambda$-Knaster. Thus the first extension does not add a new branch to $T$.

(2) $\mathbb{P}^\lambda$ is $\kappa^+$-Knaster in $W[M][j(Q)]$ where it is still Cohen forcing, and thus at least $\kappa^+$-cc in $W[M][j(Q)][\overline{P}]$ because $\mathbb{\overline{P}}$ is $\kappa^+$-cc. Since $\kappa^+<\lambda$, the second extension does not add branch.

(3) By elementarity $j(\mathbb{\overline{P}})/\overline{P}$ is $\kappa^+$-cc in $W[j(M)\times j(Q)][\overline{P}]$, and thus in the smaller model $W[M][P^\lambda][j(Q)][\overline{P}]$, so the third extension does not add branch.

(4) $\mathbb{R}^\lambda$ is $\lambda$-closed in $W[M]$. Applying Lemma \ref{forcing_branch_Unger_formerly_closed} in $W[M]$ to $(\mathbb{P}^\lambda\times j(\mathbb{Q})) *j(\mathbb{\overline{P}})$ and $\mathbb{R}^\lambda$, we see that the fourth extension does not add branch. Note that technically $j(\mathbb{\overline{P}})$ is a $j(\mathbb{P})$-name in $V$ but we are viewing it as a $\mathbb{P}^\lambda$-name in $V[P]$.
\end{proof}

\section{Questions}

In view of Theorem \ref{KH_comb}, an obvious question is whether all combinations of $\{\mathsf{KH}(\aleph_m,\aleph_n):1\leq m\leq n\le 3\}$ that do not violate Lemma \ref{KH_constraint} are consistent (there are $40$ of them). There seem to be a large number of cases for which our current method is insufficient. For example, we do not know the answer to the following

\begin{question}
Is $-\mathsf{KH}(\aleph_1,\aleph_3)-\mathsf{KH}(\aleph_2,\aleph_2)+\mathsf{KH}(\aleph_2,\aleph_3)$ consistent?
\end{question}

This is related to a limitation of our proof of Theorem \ref{KH_comb} (ii), which forces over $L$. It is natural to ask

\begin{question}
Assuming suitable large cardinals, is it always possible to force $-\mathsf{KH}(\aleph_1,\aleph_1)+\mathsf{KH}(\aleph_1,\aleph_2)+\mathsf{KH}(\aleph_2,\aleph_2)$?\footnote{We have a tentative positive answer and are working on the details.}
\end{question}

Our proof of Theorem \ref{KH_comb} (v) uses $\mathbb{Q}(\aleph_1,\aleph_2,\kappa)$, which by Lemma \ref{KH_Q_effect} blows up the power of $\aleph_1$ to $\kappa$, and hence $2^{\aleph_1}=\aleph_3$ in the final model.

\begin{question}
Is $\mathsf{GCH}+\mathsf{KH}(\aleph_1,\aleph_1)+\mathsf{KH}(\aleph_1,\aleph_2)-\mathsf{KH}(\aleph_2,\aleph_2)$ consistent?
\end{question}

Another limitation of the forcing $\mathbb{Q}(\kappa^+,\lambda,\mu)$ is the assumption $\kappa^{<\kappa}=\kappa$, which is needed for canonical $\kappa$-closure. We do not know if there is an appropriate version of the forcing that works for singular $\kappa$. Below is a test question.

\begin{question}
Is $-\mathsf{KH}(\aleph_1,\kappa^+)+\mathsf{KH}(\kappa^+,\kappa^+)$ consistent for $\kappa$ singular?
\end{question}

\bibliographystyle{plain}
\bibliography{ref}

\begin{thebibliography}{10}

\bibitem{suslin1920}
Probl\`{e}mes.
\newblock {\em Fundamenta Mathematicae}, 1:223--224, 1920.

\bibitem{cody-white2024}
Brent Cody and Philip White.
\newblock Two-cardinal ideal operators and indescribability.
\newblock {\em Annals of Pure and Applied Logic}, 175(8):103463, 2024.

\bibitem{cummings_notes}
James Cummings.
\newblock Upper bounds.
\newblock unpublished notes.

\bibitem{cummings2018aronszajn}
James Cummings.
\newblock Aronszajn and kurepa trees.
\newblock {\em Archive for Mathematical Logic}, 57(1):83--90, 2018.

\bibitem{cummings1998tree}
James Cummings and Matthew Foreman.
\newblock The tree property.
\newblock {\em Advances in Mathematics}, 133(1):1--32, 1998.

\bibitem{devlin1973kurepa}
Keith~J Devlin.
\newblock Kurepa's hypothesis and the continuum.
\newblock {\em Preprint series: Pure mathematics http://urn. nb. no/URN: NBN:
  no-8076}, 1973.

\bibitem{devlin2017constructibility}
Keith~J Devlin.
\newblock {\em Constructibility}, volume~6.
\newblock Cambridge University Press, 2017.

\bibitem{donder-matet1993}
Hans-Dieter Donder and Pierre Matet.
\newblock Two cardinal versions of diamond.
\newblock {\em Israel Journal of Mathematics}, 83(1):1--43, 1993.

\bibitem{friedman-golshani2012independence}
Sy-David Friedman and Mohammad Golshani.
\newblock Independence of higher kurepa hypotheses.
\newblock {\em Archive for Mathematical Logic}, 51(5):621--633, 2012.

\bibitem{golshani2019generalized}
Mohammad Golshani.
\newblock The generalized kurepa hypothesis at singular cardinals.
\newblock {\em Periodica Mathematica Hungarica}, 78(2):200--202, 2019.

\bibitem{hayut2017chang}
Yair Hayut.
\newblock Magidor--malitz reflection.
\newblock {\em Archive for Mathematical Logic}, 56(3):253--272, 2017.

\bibitem{honzik2020indestructibility}
Radek Honzik and {\v{S}}{\'a}rka Stejskalov{\'a}.
\newblock Indestructibility of the tree property.
\newblock {\em The Journal of Symbolic Logic}, 85(1):467--485, 2020.

\bibitem{jech1973mess}
Thomas~J Jech.
\newblock Some combinatorial problems concerning uncountable cardinals.
\newblock {\em Annals of Mathematical Logic}, 5(3):165--198, 1973.

\bibitem{Jensen-Kunen}
Ronald~B Jensen and Kenneth Kunen.
\newblock Some combinatorial properties of {L} and {V}.
\newblock handwritten notes.

\bibitem{Kanamori_2011}
Akihiro Kanamori.
\newblock {\em Historical remarks on Suslin’s problem}, page 1–12.
\newblock Lecture Notes in Logic. Cambridge University Press, 2011.

\bibitem{kurepa1936}
Djuro~R. Kurepa.
\newblock L’hypoth\`{e}se de ramification.
\newblock {\em Comptes rendues hebdomadaires des s\'{e}ances de
  l’Acad\'{e}mie des Sciences, Paris}, 202:185--187, 1936.

\bibitem{levinski1990chang}
Jean-Pierre Levinski, Menachem Magidor, and Saharon Shelah.
\newblock Chang’s conjecture for $\aleph_\omega$.
\newblock {\em Israel Journal of Mathematics}, 69(2):161--172, 1990.

\bibitem{magidor1974}
Menachem Magidor.
\newblock Combinatorial characterization of supercompact cardinals.
\newblock {\em Proceedings of the American Mathematical Society}, pages
  279--285, 1974.

\bibitem{sharon-viale2010}
Assaf Sharon and Matteo Viale.
\newblock Some consequences of reflection on the approachability ideal.
\newblock {\em Transactions of the American Mathematical Society},
  362(8):4201--4212, 2010.

\bibitem{shelah-shioya2006}
Saharon Shelah and Masahiro Shioya.
\newblock Nonreflecting stationary sets in $\mathcal{P}_\kappa\lambda$.
\newblock {\em Advances in Mathematics}, 199(1):185--191, 2006.

\bibitem{todorcevic2007walks}
Stevo Todorcevic.
\newblock {\em Walks on ordinals and their characteristics}.
\newblock Springer, 2007.

\bibitem{unger2012fragility}
Spencer Unger.
\newblock Fragility and indestructibility of the tree property.
\newblock {\em Archive for Mathematical Logic}, 51(5):635--645, 2012.

\bibitem{unger2015fragility}
Spencer Unger.
\newblock Fragility and indestructibility ii.
\newblock {\em Annals of Pure and Applied Logic}, 166(11):1110--1122, 2015.

\bibitem{weiss2010subtle}
Christoph Wei{\ss}.
\newblock {\em Subtle and ineffable tree properties}.
\newblock PhD thesis, lmu, 2010.

\end{thebibliography}

\end{document}